\documentclass[12pt]{article}

\usepackage[utf8]{inputenc}
\usepackage[margin=3cm]{geometry}
\usepackage{amsfonts, amsmath, amssymb, amsthm}
\usepackage{indentfirst}
\usepackage{xcolor}

\newtheorem{Thm}{Theorem}
\newtheorem{Prp}{Proposition}
\newtheorem{Lem}{Lemma}

\newenvironment{Prf}{\par\noindent{\it Proof:}\rm}{\newline\rightline{\textbf{QED}}\par}
\newtheorem{Cor}{Corollary}

\newcommand{\R}{\mathbb{R}}

\newcommand{\con}{\mathfrak{c}}

\newcommand{\Perf}{\rm{Perf}}
\newcommand{\Pow}{\mathcal{P}}

\newcommand{\dom}{\mathrm{dom}}

\title{
  \textbf{On the methods of reduction of some types of Marczerwski-Burstin measurable functions to continuous functions on products of perfect sets}
}

\author{Waldemar Hołubowski, Sławomir Kusiński}

\date{}

\begin{document}
\maketitle

\begin{abstract}
  In this paper, we introduce product-wise generalizations of certain Marczewski-Burstin bases, including sets with the (s)-property and completely Ramsey sets.
  For each of these families, we establish analogs of the classical Luzin and Eggleston theorems,
  showing that functions measurable with respect to these families can be reduced to continuous functions on products of perfect sets.
  Furthermore, we provide a method for reducilng sequences of such functions to continuity,
  which allows us to generalize Laver's extension of Halpern-Läuchli and Harrington theorems.
\end{abstract}


%%%%%%%%%%%%%%%%%%%%
%
% INTRODUCTION
%
%%%%%%%%%%%%%%%%%%%%
\section{Introduction}

In~\cite{Bur}, Burstin showed that if we consider the families
\begin{equation}\label{MarBur}
  \mathcal{S}(\mathcal{F}) = \{ S \subseteq X \colon \forall_{P \in \mathcal{F}}\exists_{Q \in \mathcal{F}} Q \subseteq S \cap P \mbox{ or } Q \subseteq (X \setminus S) \cap P \}
\end{equation}
and
\begin{equation}
  \mathcal{S}_0(\mathcal{F}) = \{ S \subseteq X \colon \forall_{P \in \mathcal{F}}\exists_{Q \in \mathcal{F}} Q \subseteq (X \setminus S) \cap P \},
\end{equation}
where \(X\) is the real line and \( \mathcal{F} \) is a family of all perfect subsets of \(X\) with positive Lebesgue measure,
then \(\mathcal{S}(\mathcal{F})\) is the family of all Lebesgue measurable sets,
and \(\mathcal{S}_0(\mathcal{F})\) is the family of all Lebesgue null sets.

In~\cite{Mar}, Marczewski introduced the notion of the \((s)\)-property, a construction closely related to Burstin's work.
Specifically, a set \(A \subseteq X \) has the \( (s) \)-property
if for any non-empty perfect set \( P \subseteq X \)
there exists a non-empty perfect set \( Q \subseteq P \) such that
\begin{equation}
  Q \subseteq A \mbox{ or } Q \subseteq X \setminus A.
\end{equation}
Similarly, \( A \) has the nowhere \((s)\)-property
if for any non-empty perfect set \( P \subseteq X \)
there exists a non-empty perfect set \( Q \subseteq P \) such that
\begin{equation}
  Q \cap A = \emptyset.
\end{equation}
The families of subsets of \(X\) with the \((s)\)-property and nowhere \((s)\)-property -- denoted after \cite{Paradise} \(s\) and \(s_0\) respectively --
are in fact \( \mathcal{S}(\Perf(X)) \) and \( \mathcal{S}_0(\Perf(X)) \), where \( \Perf(X) \) denotes the family of all perfect subsets of \(X\).
Just like Burstin, Marczewski considered the case where \( X = \R \), but \(X\) can be taken to be any dense in itself Polish space.
The construction~\ref{MarBur} in general is often called Marczewski-Burstin base~\cite{MarBur}.
Clearly, the family \( \mathcal{S}(\mathcal{F}) \) is always an algebra of sets, and \( \mathcal{S}_0(\mathcal{F}) \) is an ideal within \( \mathcal{S}(\mathcal{F}) \).

Another well-known example of Marczewski-Burstin base is the notion of completely Ramsey and nowhere Ramsey sets.
In~\cite{El}, Ellentuck has introduced a very interesting topology on \( [\omega]^\omega \) in a following way.
The base of the topology consists of the sets
\begin{equation*}
  (s, A)^\omega = \{ B \in [\omega]^\omega \colon s \subseteq B \subseteq A \cup s, \max(s) < \min(B \setminus s)  \},
\end{equation*}
where \( A \in [\omega]^\omega \) and \( s \in [A]^{<\omega} \) with \( \max(s) < \min(A) \).
It is easy to see that such a topology is stronger than the euclidean topology, i.e.~the subspace topology when considering \( [\omega]^\omega \) as a subset of \( 2^\omega \).
It has been shown in~\cite{Plew1} that such a topological space is not metrizable and furthermore it is not even normal.
The families of completely and nowhere Ramsey sets are \( r = \mathcal{S}(\mathcal{U}) \) and \( r_0 = \mathcal{S}_0(\mathcal{U}) \),
where \(\mathcal{U}\) is the family of all the sets \( (s, A)^\omega \).
The nowhere Ramsey sets coincide with sets that are nowhere dense in the Ellentuck topology.

The function \( f \colon X \to Y \) where \(Y\) is a metric space is called measurable with respect to the family \(\mathcal{S}(\mathcal{F})\)
if for any open set \( U \subseteq Y \) we have \( f^{-1}(U) \in \mathcal{S}(\mathcal{F}) \).
In~\cite{RFAni} the authors proved the following two theorems, which could be viewed as generalizations of Luzin's theorem~\cite{Luz}.
\begin{Thm}
  A function \( f \colon \R \to Y \) is \( (s) \)-measurable iff for any perfect set \( P \subseteq \R \) there exists a perfect set \( Q \subseteq P \) such that \( f|_Q \) is continuous.
\end{Thm}
\begin{Thm}
  A function \( f \colon [\omega]^{\omega} \to Y \) is completely Ramsey (i.e. \( (r) \)-measurable) iff for any base set \( (s, A)^\omega \) there exists \( B \in [A]^{\omega} \) such that \( f|_{(s, B)^\omega} \) is continuous with respect to the euclidean topology.
\end{Thm}

In this paper we generalize these results to three particularly interesting classes of product spaces, in a manner analogous to how Eggleston~\cite{Egg} generalized Luzin's theorem. Our main theorems will be the following.
\begin{Thm}\label{thm:s-property}
  Let \( X_m \) be dense-in-itself Polish spaces for \( m \in d \), where \( d \le \omega \).
  If \( f \colon \prod\limits_{m \in d} X_m \to Y \) is \( (s^d) \)-measurable, then for any perfect \(d\)-cube \( P = \prod\limits_{m \in d} P_m \)
  there exists a perfect \(d\)-cube \( Q \subseteq P \), such that \( f|_Q \) is continuous.
\end{Thm}
\begin{Thm}\label{thm:v-property}
  Let \( d \le \omega \).
  If \( f \colon (2^\omega)^d \to Y \) is \( (v^d) \)-measurable, then for any Silver \(d\)-cube \( P = \prod\limits_{m \in d} [h_m] \)
  there exists a Silver \(d\)-cube \( Q \subseteq P \), such that \( f|_Q \) is continuous.
\end{Thm}
\begin{Thm}\label{thm:r-property}
  Let \( d \le \omega \).
  If \( f \colon ([\omega]^\omega)^d \to Y \) is \( (r^d) \)-measurable, then for any Ellentuck \(d\)-cube \( P = \prod\limits_{m \in d} (s_m, A_m)^\omega \)
  there exists an Ellentuck \(d\)-cube \( Q \subseteq P \), such that \( f|_Q \) is continuous with respect to the euclidean topology.
\end{Thm}

Beyond proving these generalizations, we also demonstrate their applicability
by extending a classical result.
Specifically, we show how these theorems enable a generalization of Laver’s
theorem~\cite{Lav} to new classes of sets, which itself was a refinement of
earlier results by Harrington and by Halpern and Läuchli~\cite{HL}.
This connection highlights the broader relevance of our approach and its
connection to a wider context of descriptive set theory
and infinite combinatorics.

%%%%%%%%%%%%%%%%%%%%
%
% PRODUCT-WISE MARCZEWSKI BURSTIN BASES
%
%%%%%%%%%%%%%%%%%%%%
\section{Product-wise Marczewski Burstin bases}

For our purposes, we will extend the notion of the \((s)\)-property and \((s)\)-measurability to accommodate products.
Let \(X_m\) be dense-in-itself Polish spaces for \( m \in d \) and \( d \le \omega \).
A set \( P = \prod\limits_{m \in d} P_m \), where all the \( P_m \subseteq X_m \) are perfect, will be called a perfect \(d\)-cube.
In the case \( d = \omega \), we will simply refer to \(P\) as a perfect cube.
Let \( A \subseteq \prod\limits_{m \in d} X_m \).
We say that \(A\) has the \( (s^d) \)-property (denoted \(A \in s^d\)),
if for any perfect \(d\)-cube \( P \), there exists a perfect \(d\)-cube \( Q \subseteq P \), such that either
\[
  Q \subseteq A \mbox{ or } Q \subseteq \prod\limits_{m \in d} X_m \setminus A.
\]
Analogously, \( A \) has the nowhere \( (s^d) \)-property (denoted \( A \in s^d_0 \)),
if for any perfect \(d\)-cube \( P \), there exists a perfect \(d\)-cube \( Q \subseteq P \), such that
\[
  Q \subseteq \prod\limits_{m \in d} X_m \setminus A.
\]

The following fact about sets with the \((s^d)\)-property will be useful to us later on.
\begin{Prp}
  The family \( s^d \) forms a \( \sigma \)-algebra, and \( s^d_0 \) is a \( \sigma \)-ideal within it.
\end{Prp}
\begin{Prf}
  Let \( A_n \in s^d \) for \( n \in \omega \), and let \( P \) be a perfect \( d \)-cube. If for some \( n \in \omega \) there exists a perfect \(d\)-cube \(Q \subseteq P \) such that \( Q \subseteq A_n \), then clearly
  \[
    Q \subseteq \bigcup\limits_{n \in \omega} A_n
  \]
  as required.

  Otherwise, for every \( n \in \omega \) and any perfect \(d\)-cube \(Q \subseteq P \), there exists a perfect \(d\)-cube \( R \subseteq Q \) such that
  \[
    R \cap \bigcup\limits_{n \in \omega} A_n = \emptyset.
  \]
  It remains to show that there exists a perfect \(d\)-cube \(Q\) such that
  \[
    Q \cap \bigcup\limits_{n \in \omega} A_n = \emptyset
  \]
  as it will simultaneously show that if all \( A_n \in s^d_0 \) then \( \bigcup\limits_{n \in \omega} A_n \in s^d_0 \).

  If \( d < \omega \) proceed as follows.
  Let \( Q^m_\emptyset = P_m \).
  Now suppose we have defined perfect sets \( Q^m_t \subseteq P_m \) for \( t \in 2^n \) and \( m \in d \) each of diameter less than \( \frac{1}{2^n} \) and satisfying
  \[
    Q^0_{t_0} \times \ldots \times Q^{d - 1}_{t_{d-1}} \cap \bigcup\limits_{k \in n} A_k = \emptyset.
  \]
  As \( A_n \in s^d_0 \) and there are only finitely many cubes  \( Q^0_{t_0} \times \ldots \times Q^{d - 1}_{t_{d-1}}\), there have to exist sets \( Q^{*m}_{t} \subseteq Q^{m}_{t} \) such that
  \[
    Q^{*0}_{t_0} \times \ldots \times Q^{*d - 1}_{t_{d-1}} \cap A_{n} = \emptyset.
  \]
  In each set \( Q^{*m}_{t} \) pick two disjoint perfect subsets \( Q^{m}_{t^\smallfrown 0}, Q^{m}_{t^\smallfrown 1} \) each of diameter less than \( \frac{1}{2^n} \).

  From the construction we obtain that the set
  \[
    Q = \prod\limits_{m \in d} \bigcap\limits_{n \in \omega} \bigcup\limits_{t \in 2^n} Q^m_{t}
  \]
  is a perfect \(d\)-cube and we have
  \[
    Q \cap \bigcup\limits_{n \in \omega} A_n = \emptyset
  \]
  as required.

  In case \( d = \omega \), let \( Q^0_\emptyset = P_0 \) and \( R^m_0 = P_m \) for \( m > 0 \).
  Now suppose we have defined perfect sets \( Q^m_t \subseteq P_m \) for \( t \in 2^n \) and \( m \le n \)
  each of diameter less than \( \frac{1}{2^n} \),
  as well as perfect sets \( R^m_n \) for \( m > n \) such that
  \[
    Q^0_{t_0} \times \ldots \times Q^{n}_{t_{n}} \times \prod\limits_{m > n} R^m_n \cap \bigcup\limits_{k \in n} A_k = \emptyset.
  \]
  Once again, as there are only finitely many cubes \( Q^0_{t_0} \times \ldots \times Q^{n}_{t_{n}} \times \prod\limits_{m > n} R^m_n \),
  there have to exist sets \( Q^{*m}_{t} \subseteq Q^{m}_{t} \) and \( R^m_{n+1} \subseteq R^m_{n} \) such that
  \[
    Q^{*0}_{t_0} \times \ldots \times Q^{*n}_{t_{n}} \times \prod\limits_{m > n} R^m_{n+1} \cap A_n = \emptyset.
  \]
  In each set \( Q^{*m}_{t} \) pick two disjoint perfect subsets \( Q^{m}_{t^\smallfrown 0}, Q^{m}_{t^\smallfrown 1} \)
  each of diameter less than \( \frac{1}{2^n} \),
  and in the set \( R^{n+1}_{n+1} \) pick \( 2^{n+1} \) many pairwise disjoint perfect subsets \( Q^{n+1}_{t} \)
  each of diameter less than \( \frac{1}{2^n} \).

  From the construction we obtain that the set
  \[
    Q = \prod\limits_{m \in \omega} \bigcap\limits_{n > m} \bigcup\limits_{t \in 2^n} Q^m_{t}
  \]
  is a perfect cube, and we have
  \[
    Q \cap \bigcup\limits_{n \in \omega} A_n = \emptyset
  \]
  as required.
\end{Prf}

We will also consider a close variant of the \((s)\)-property based on the Prikry-Silver~forcing~\cite{JechM}.
Let \( D \subseteq \omega \) and \( h \colon D \to 2 \).
Then the subset of the Cantor set \( 2^\omega \) generated by \(h\) is defined as \( [h] = \{ x \in 2^\omega \colon \forall_{n \in D} x(n) = h(n) \} \).
A set \( P = \prod\limits_{m \in d} [h_m] \subseteq (2^\omega)^d \),
where each \( h_m \in 2^{D_m}  \) and \( D_m \) is a coinfinite subset of \( \omega \),
is called a Silver \(d\)-cube.
When \( d = \omega \), we refer to \(P\) simply as a Silver cube.
Clearly, each Silver \(d\)-cube is a perfect \(d\)-cube as well.
Similarly to the \((s)\)-property, we say that \(A \subseteq (2^\omega)^d \) has the \( (v^d) \)-property (denoted \(A \in v^d\))
if for any Silver \(d\)-cube \( P \), there exists a Silver \(d\)-cube \( Q \subseteq P \) such that either
\[
  Q \subseteq A \mbox{ or } Q \subseteq (2^\omega)^d \setminus A.
\]
Analogously, \( A \) has the nowhere \( (v^d) \)-property (denoted \( A \in v^d_0 \)) if for any Silver \(d\)-cube \( P \), there exists a Silver \(d\)-cube \( Q \subseteq P \) such that
\[
  Q \subseteq (2^\omega)^d \setminus A.
\]

\vspace{20pt}

The \( (v^d) \)-property turns out to be \(\sigma\)-additive as well.
\begin{Prp}
  The family \( v^d \) forms a \( \sigma \)-algebra, and \( v^d_0 \) is a \( \sigma \)-ideal within it.
\end{Prp}
\begin{Prf}
  Let \( A_n \in v^d \) for \( n \in \omega \), and let \( P \) be a Silver \(d\)-cube.
  If for some \( n \in \omega \) there exists a Silver \(d\)-cube \(Q \subseteq P \) such that \( Q \subseteq A_n \), then clearly
  \[
    Q \subseteq \bigcup\limits_{n \in \omega} A_n
  \]
  as required.

  It remains to show that if for any \( n \in \omega \) and any Silver \(d\)-cube \(Q \subseteq P \) there exists a Silver \(d\)-cube \(R \subseteq Q \) such that
  \[
    R \cap A_n = \emptyset,
  \]
  then there exists a Silver \(d\)-cube \( Q \subseteq P \) such that
  \[
    Q \cap \bigcup\limits_{n \in \omega} A_n = \emptyset.
  \]
  We will need to consider two separate cases.

  \noindent(\(d < \omega\))
  \\
  Let \( Q_0 = \prod\limits_{m \in d} [h_{m,0}] = P \), and let \( i_{m,0} = \min (\omega \setminus \dom(h_{m,0})) \).
  Now suppose we have defined a Silver \(d\)-cube \( Q_n = \prod\limits_{m \in d} [h_{m,n}] \subseteq P \)
  and numbers \( i_{m,0}, \ldots, i_{m,n} \not\in \dom(h_{m,n}) = D_{m,n} \) such that
  \[
    Q_n \cap \bigcup\limits_{k \in n} A_k = \emptyset.
  \]
  Let \( (2^n)^m = \{ \alpha_k \colon k \in 2^{n \cdot m} \} \) and \( Q^{0}_{m,n} = [h^{0}_{m, n}] = Q_{m,n} \).
  For any \( k \in 2^n \) let
  \[
    h^{*k+1}_{m, n} = h^k_{m, n}|_{\omega \setminus I_{m,n})} \cup \{ (i_{m,0}, \alpha_{k}(m,0)), \ldots, (i_{m,n}, \alpha_{k}(m,n)) \},
  \]
  where \( I_{m,n} = \{ i_{m,0}, \ldots, i_{m,n} \} \).
  In a Silver \(d\)-cube \( \prod\limits_{m \in d} [h^{*k+1}_{m, n}] \) pick a Silver \(d\)-cube \( \prod\limits_{m \in d} [h^{k+1}_{m, n}] \) disjoint from \( A_n \).
  We can take \( Q_{m, n+1} = \prod\limits_{m \in d} [h_{m, n+1}] = [h^{2^{n \cdot d}}_{m, n}|_{\omega \setminus I_{m,n}}] \) and \( i_{m,n+1} = \min (\omega \setminus D_{m,n+1}) \)
  where \( D_{m,n+1} = \dom (h_{m,n+1}) \).
  Clearly \( i_{m,n+1} > i_{m,n} \).
  As each set \( I_m = \bigcup\limits_{n \in \omega} I_{m,n} \) is infinite and disjoint from any \( \dom(h_{m,n+1}) \), we obtain that the set
  \[
    Q = \prod\limits_{m \in d} \bigcap\limits_{n \in \omega} Q_{m,n} = \prod\limits_{m \in d} [\bigcup\limits_{n \in \omega} h_{m,n}]
  \]
  is a Silver \(d\)-cube, and it is disjoint from any \( A_n \) as required.

  \noindent(\(d = \omega\))
  \\
  Let \( Q_0 = \prod\limits_{m \in \omega} [h_{m,0}] = P \), and let \( i_{0,0} = \min (\omega \setminus \dom(h_{0,0})) \).
  Now suppose we have defined a Silver cube \( Q_n = \prod\limits_{m \in \omega} [h_{m,n}] \subseteq P \)
  and numbers \( i_{m,0}, \ldots, i_{m,n} \not\in \dom(h_{m,n}) = D_{m,n} \) for \( m \le n \) such that
  \[
    Q_n \cap \bigcup\limits_{k \in n} A_k = \emptyset.
  \]
  Let \( (2^n)^n = \{ \alpha_k \colon k \in 2^{n^{2}} \} \) and \( Q^{0}_{m,n} = [h^{0}_{m, n}] = Q_{m,n} \).
  For any \( k \in 2^n \) let
  \[
    h^{*k+1}_{m, n} = h^k_{m, n}|_{(\omega \setminus I_{m,n})}
    \cup \{ (i_{m,0}, \alpha_{k}(m,0)), \ldots, (i_{m,n}, \alpha_{k}(m,n)) \}
  \]
  for \( m \le n\), and \(h^{*k+1}_{m, n} = h^{k}_{m, n}\) for \( m > n \),
  where \( I_{m,n} = \{ i_{m,0}, \ldots, i_{m,n} \} \).
  In a Silver cube \( \prod\limits_{m \in \omega} [h^{*k+1}_{m, n}] \) pick a Silver cube \( \prod\limits_{m \in \omega} [h^{k+1}_{m, n}] \) disjoint from \( A_n \).
  We can take \( Q_{m, n+1} = \prod\limits_{m \in \omega} [h_{m, n+1}] = [h^{2^{n \cdot m}}_{m, n}|_{\omega \setminus I_{m,n}}] \) and \( i_{m,n+1} = \min (\omega \setminus D_{m,n+1}) \)
  where \( D_{m,n+1} = \dom (h_{m,n+1}) \).
  Clearly \( i_{m,n+1} > i_{m,n} \).
  Furthermore, take \( i_{0,n+1}, \ldots, i_{n+1, n+1} \) to be the first \( n+1 \) natural numbers not in \( D_{n+1,n+1} \).
  As each set \( I_m = \bigcup\limits_{n \ge m} I_{m,n} \) is infinite and disjoint from any \( \dom(h_{m,n+1}) \), we obtain that the set
  \[
    Q = \prod\limits_{m \in \omega} \bigcap\limits_{n \in \omega} Q_{m,n} = \prod\limits_{m \in \omega} [\bigcup\limits_{n \in \omega} h_{m,n}]
  \]
  is a Silver cube, and it is disjoint from any \( A_n \) as required.
\end{Prf}

The last product-wise generalization of Marczewski-Burstin base we consider involves completely Ramsey sets.
We define a set of the form \( \prod\limits_{n \in d} (s_n, A_n)^\omega \) as an Ellentuck \(d\)-cube.
A set \( S \subseteq ([\omega]^\omega)^d \) has the \( (r^d) \)-property if for any Ellentuck \(d\)-cube \( U \), there exists an Ellentuck \(d\)-cube \( V \subseteq U \) such that either
\[
  S \subseteq V \mbox{ or } S \subseteq ([\omega]^\omega)^d \setminus V,
\]
and it has the \( (r^d_0) \)-property if for any Ellentuck \(d\)-cube \( U \) there exists an Ellentuck \(d\)-cube \( V \subseteq U \) such that
\[
  S \subseteq ([\omega]^\omega)^d \setminus V.
\]

There is one simple property of Ellentuck base sets that will be useful to us later on.
\begin{Prp}
  For any \( s \in [\omega]^{<\omega} \) and \( A \in [\omega]^{\omega} \) such that \( \max(s) < \min(A) \) we have
  \[
    \bigcup\limits_{\alpha \in \Pow(s)} (\alpha, A)^\omega = (\emptyset, A \cup s )^\omega
  \]
\end{Prp}
\begin{Prf}
  Clearly, \( \bigcup\limits_{\alpha \in \Pow(s)} (\alpha, A)^\omega \subseteq (\emptyset, A \cup s )^\omega \), as \( (\emptyset, A \cup s )^\omega = [A \cup s]^\omega \).
  Let \( B \in [A \cup s]^\omega \), and define \( \alpha = B \cap s \). 
  Then \( \max (\alpha) < \min(A) \), and since \( B \subseteq \alpha \cup A \), it follows that \( B \in  (\alpha, A)^\omega \).
\end{Prf}

For our purposes, \(\sigma\)-additivity of the ideal \( r^d_0 \) will be needed.
\begin{Prp}
  Let \( S_n \in r^d_0 \) for \( n \in \omega \). Then \( \bigcup\limits_{n \in \omega} S_n \in r^d_0 \).
\end{Prp}
\begin{Prf}
  Let \( \prod\limits_{m \in d} (s_m, A_m)^\omega \) be any Ellentuck \(d\)-cube
  and \( A_{m,0} = A_{m}\).
  Without loss of generality we can assume \( s_m = \emptyset \).

  First consider the case \( d \in \omega \).
  Suppose we have defined the sets \( A_{m,n} \) as~well as~numbers~\( a_{m,k} \) for \( m \in d \) and \( k \in n \), such that
  \( \max(\{ a_{m,0}, \ldots, a_{m, n-1} \}) < \min(A_{m,n}) \) and
  \[
    \prod\limits_{m \in d} (\emptyset, A_{m,n} \cup \{ a_{m,0}, \ldots, a_{m, n-1} \})^\omega \cap S_k = \emptyset
  \]
  for \( k \in n \).
  Let \( \prod\limits_{m \in d}\Pow(\{ a_{m,0}, \ldots, a_{m, n-1} \}) = \{ \alpha^k \colon k \in 2^{d \cdot n} \} \) and \( A^0_{m,n} = A_{m,n} \).\\
  For~any~\(k \in d \cdot n\) in the Ellentuck \(d\)-cube \(\prod\limits_{m \in d} (\alpha^k_m, A^k_{m,n})^\omega\)
  find  the Ellentuck \(d\)-cube \(\prod\limits_{m \in d} (\alpha^k_m, A^{k+1}_{m,n})^\omega\)
  disjoint with \(S_n\).
  Let \( a_{m,n} = \min A^{d \cdot n}_{m,n} \).
  We can take \( A_{m,n+1} = A^{2^{d \cdot n}}_{m,n} \setminus \{ a_{m,n} \} \).

  From the contruction we obtain that the Ellentuck \(d\)-cube
  \[
    \prod\limits_{m \in d} (\emptyset, \{ a_{m, n} \colon n \in \omega \})^\omega
  \]
  is disjoint with every \( S_k \).

  In case when \( d = \omega \), suppose we have defined the sets \( A_{m,n} \) for \( m \in \omega \) and numbers \( a_{m,k} \) for \( m, k \in n \), such that
  \( \max(\{ a_{m,0}, \ldots, a_{m, n-1} \}) < \min(A_{m,n}) \) and
  \[
    \prod\limits_{m \in n} (\emptyset, A_{m,n} \cup \{ a_{m,0}, \ldots, a_{m, n-1} \})^\omega \times \prod\limits_{m > n} (\emptyset, A_{m,n})^\omega \cap S_k = \emptyset
  \]
  for \( k \in n \).
  Let \( \prod\limits_{m \in n}\Pow(\{ a_{m,0}, \ldots, a_{m, n-1} \}) = \{ \alpha^k \colon k \in 2^{n^2} \} \) and \( A^0_{m,n} = A_{m,n} \).
  For any \(k \in n^2\) in the Ellentuck cube \(\prod\limits_{m \in n} (\alpha^k_m, A^k_{m,n})^\omega \times \prod\limits_{m > n} (\emptyset, A^k_{m,n})^\omega \)
  find  the Ellentuck cube \(\prod\limits_{m \in n} (\alpha^k_m, A^{k+1}_{m,n})^\omega \times \prod\limits_{m > n} (\emptyset, A^{k+1}_{m,n})^\omega \)
  disjoint with \(S_n\).
  Let \( a_{n, 0}, \ldots, a_{n,n} \) be the first \( n + 1 \) elements of \( A^{n^2}_{n,n} \),
  and for \( m \in n \) let \( a_{m,n} = \min A^{n^2}_{m,n} \).
  We can take \( A_{m,n+1} = A^{2^{n^2}}_{m,n} \setminus \{ a_{m,n} \} \) for \( m \in n \),
  \( A_{n, n+1} = A^{2^{n^2}}_{n,n} \setminus \{ a_{n, 0}, \ldots, a_{n,n} \} \),
  and \( A_{m, n+1} = A^{2^{n^2}}_{m,n} \) for \( m > n \).

  From the construction we obtain that the Ellentuck cube
  \[
    \prod\limits_{m \in \omega} (\emptyset, \{ a_{m, n} \colon n \in \omega \})^\omega
  \]
  is disjoint with every \( S_k \).
\end{Prf}
The \(\sigma\)-completeness of the algebra \(r^d\) can be proven using methods of Galvin and Prikry~\cite{GP}.

%%%%%%%%%%%%%%%%%%%%
%
% MAIN RESULTS
%
%%%%%%%%%%%%%%%%%%%%
\section{Main result}

\begin{Lem}
  Let \( \mathcal{A} = \{ A_i \colon i \in I \} \) be a disjoint,  \( (s^d) \)-additive family.
  If there exist a perfect \(d\)-cube \( P \prod\limits_{m \in d} P_m \subseteq \bigcup \mathcal{A} \)
  then the set \(\{i \in I \colon A_i \cap P \not= \emptyset \}\) has cardinality \(\con \).
\end{Lem}
\begin{Prf}
  As \(s^d_0\) is a \(\sigma\)-ideal in \( s^d \) the set \( \Delta = \{ J \subseteq I \colon \bigcup\limits_{i \in J} A_i \in s^d_0 \} \) is a free \(\sigma\)-ideal on \( I \).
  Let \(J = \{i \in I \colon A_i \cap P \not= \emptyset \}\).
  Clearly \( J \not\in \Delta \) and there exist \( J_0, J_1 \subseteq J \) disjoint and both not in \( \Delta \).
  There exist perfect \(d\)-cubes \( P_0 \subseteq \bigcup\limits_{i \in J_0} A_i \) and \( P_1 \subseteq \bigcup\limits_{i \in J_1} A_i \).

  Now assume that we have defined disjoint sets \( J_t \not\in \Delta \) for \( t \in 2^{n+1} \) as well as perfect \(d\)-cubes \( P_t \subseteq \bigcup\limits_{i \in J_t} A_i \).
  In each \( J_t \) we can find disjoint subsets \( J_{t^\smallfrown 0}, J_{t^\smallfrown 1} \)
  such that there exist perfect \(d\)-cubes \( P_{t^\smallfrown 0} \subseteq J_{t^\smallfrown 0} \cap P_t \) and \( P_{t^\smallfrown 1} \subseteq J_{t^\smallfrown 1} \cap P_t \).

  We obtain that for any \( x \in 2^\omega \) the intersection \( \bigcap\limits_{n \in \omega} P_{x|_n} \not= \emptyset \)
  and consequently \( \bigcap\limits_{n \in \omega} J_{x|_n} \not= \emptyset \).
  Thus, \( |J| = \con \).
\end{Prf}

Furthermore, as every Silver \(d\)-cube is a perfect \(d\)-cube, the proof above is also valid for \( v^d \) and \( v^d_0 \) in place of \( s^d \) and \( s^d_0 \).

\begin{Cor}
  Let \( \mathcal{A} = \{ A_i \colon i \in I \} \) be a disjoint,  \( (v^d) \)-additive family.
  If there exist a Silver \(d\)-cube \( P \prod\limits_{m \in d} P_m \subseteq \bigcup \mathcal{A} \)
  then the set \(\{i \in I \colon A_i \cap P \not= \emptyset \}\) has cardinality \(\con \).
\end{Cor}

By using a variant of the Bernstein construction, we obtain the following.

\begin{Cor}
  Let \( \mathcal{A} = \{ A_i \colon i \in I \} \subseteq s^d_0 \) be a disjoint,  \( (s^d) \)-additive family.
  Then \( \bigcup\mathcal{A} \in s^d_0 \).
\end{Cor}
\begin{Prf}
  Let \( P = \prod\limits_{m \le d} P_m \subseteq \bigcup\mathcal{A} \) be a product of perfect sets contained in \( \bigcup \mathcal{A} \),
  and let \(\{ Q_\alpha = \prod\limits_{m \le d} Q_{\alpha,m} \colon \alpha \in 2^\omega \}\) be a family of all products of perfect sets contained in \( P \).
  Pick distinct \( i_0, j_0 \in I \), so that both \( A_{i_0} \) and \( A_{j_0} \) have a nonempty intersection with \( Q_0 \).
  With \( Q_\beta, i_\beta \) and \( j_\beta \) defined for \( \beta < \alpha < 2^\omega \),
  pick distinct \( i_\alpha, j_\alpha \in I \setminus (\{ i_\beta \colon \beta < \alpha \} \cup \{ j_\beta \colon \beta < \alpha \}) \),
  such that both \( A_{i_\alpha} \) and \( A_{j_\alpha} \) have non-empty intersection with \( Q_\alpha \).
  We get that both sets \( \bigcup\{ A_{i_\alpha}\colon \alpha \in 2^\omega \}\) and \(\bigcup\{ A_{j_\alpha}\colon \alpha \in 2^\omega \} \)
  have non-empty intersection with every set \( Q_\alpha \), which contradicts the additivity of the family \( \mathcal{A} \).
\end{Prf}

Once again, the same reasoning works for Silver \(d\)-cubes.

\begin{Cor}
  Let \( \mathcal{A} = \{ A_i \colon i \in I \} \subseteq v^d_0 \) be a disjoint,  \( (v^d) \)-additive family.
  Then \( \bigcup\mathcal{A} \in v^d_0 \).
\end{Cor}

We can now proceed with proving theorems~\ref{thm:s-property} and~\ref{thm:v-property}.
\begin{Prf} of Theorem~\ref{thm:s-property}.
  Let \( f \) be \((s^d)\)-measurable, and let \( P = \prod\limits_{m \le d} P_m \subseteq \bigcup\mathcal{A} \) be a perfect \(d\)-cube.

  Consider first the case \noindent(\(d < \omega\)).
  Let \( Q_{m, \emptyset} = P_m \).
  Suppose now that we have defined the perfect sets \( Q_{m, t} \) for \( t \in 2^n \), all of diameter less than \( \frac{1}{n} \).
  Fix a cover \( \mathcal{U}_n \) of \( X \) consisting of open sets of diameter less than \( \frac{1}{n+1} \).
  Since every metric space is paracompact, the cover \( \mathcal{U}_n \) has a \(\sigma\)-discrete refinement \( \tilde{\mathcal{U}}_n \).
  By the corollary above and \(\sigma\)-discreteness of \( \tilde{\mathcal{U}}_n \),
  for any perfect \(d\)-cube \( R \) there has to exist \( U \in \tilde{\mathcal{U}}_n \), such that \( U \cap R \in s^d \setminus s^d_0 \).
  As there are finitely many products of the form \( \prod\limits_{m < d} Q_{m, t_m} \),
  we obtain the perfect sets \( Q^*_{m, t} \subseteq Q_{m, t}  \),
  such that each product is contained in a set \( f^{-1}(U) \) for some \( U \in \tilde{\mathcal{U}} \).
  In each \( Q^*_{m,t} \) we can find two disjoint perfect subsets \( Q_{m,t^\smallfrown 0}, Q_{m,t^\smallfrown 1} \) each having diameter less than \( \frac{1}{n+1} \).

  Hence, the set
  \[
    Q = \prod\limits_{m \in d} \bigcap\limits_{n \in \omega} \bigcup\limits_{t \in 2^n} Q_{m,t}
  \]
  is a perfect \(d\)-cube.
  Furthermore, it is clear from the construction that the sets \( Q \cap \prod\limits_{m < d} Q_{m, t_m} \)
  form a base for the product topology on \(Q\).
  It follows that the function \( f|_Q \) is continuous.

  For the case \((d=\omega)\) let \(R_{m,0} = P_m\).
  Suppose now that we have defined the perfect sets \( Q_{m, t} \) for \( t \in 2^n \) and \( m \in n \), all of the diameter less than \( \frac{1}{n} \),
  as well as perfect sets \( R_{m,n} \) for \(m \ge n\).
  Exactly as in finitely dimensional case
  fix a cover \( \mathcal{U}_n \) of \( X \) consisting of open sets of diameter less than \( \frac{1}{n+1} \),
  and take its \(\sigma\)-discrete refinement \( \tilde{\mathcal{U}}_n \).
  By the corollary above and \(\sigma\)-discreteness of \( \tilde{\mathcal{U}}_n \),
  for any perfect cube \( S \) there has to exist \( U \in \tilde{\mathcal{U}}_n \), such that \( U \cap S \in s^\omega \setminus s^\omega_0 \).
  As there are finitely many products of the form \( \prod\limits_{m < n} Q_{m, t_m} \times \prod\limits_{m \ge n} R_{m,n} \),
  we obtain the perfect sets \( Q^*_{m, t} \subseteq Q_{m, t}  \) and \( R_{m, n+1} \subseteq R_{m,n} \),
  such that each product \( \prod\limits_{m < n} Q^*_{m, t_m} \times \prod\limits_{m \ge n} R_{m,n+1} \)
  is contained in a set \( f^{-1}(U) \) for some \( U \in \tilde{\mathcal{U}} \).
  In each \( Q^*_{m,t} \) we can find two disjoint perfect subsets \( Q_{m,t^\smallfrown 0}, Q_{m,t^\smallfrown 1} \) each having diameter less than \( \frac{1}{n+1} \),
  and in \(R_{n,n+1}\) we can find \( 2^{n+1} \) many disjoint perfect subsets \( Q_{n, t} \) each having diameter less than \( \frac{1}{n+1} \).

  As a result we obtain the perfect cube
  \[
    Q = \prod\limits_{m \in d} \bigcap\limits_{n \ge m} \bigcup\limits_{t \in 2^n} Q_{m,t} = \prod\limits_{m \in d} Q_m.
  \]
  Moreover, it is clear from the construction that the sets
  \[
    Q \cap \left( \prod\limits_{m < n} Q_{m, t_m} \times \prod\limits_{m \ge n} R_{m,n} \right)
    = Q \cap \left( \prod\limits_{m < n} Q_{m, t_m} \times \prod\limits_{m \ge n} Q_{m} \right)
  \]
  form a base for the product topology on \(Q\).
  It follows that the function \( f|_Q \) is continuous.
\end{Prf}

\begin{Prf} of Theorem \ref{thm:v-property}.
  Let \( f \) be \((v^d)\)-measurable and \( P = \prod\limits_{m \le d} P_m = \prod\limits_{m \le d} [h_m] \subseteq \bigcup\mathcal{A} \) be a Silver \(d\)-cube.

  First consider the case \(d < \omega\).
  Let \(Q_{m,0} = P_m\) and \(i_{m,0} = \min(\omega \setminus \dom(h_m)) \).
  Suppose now that we have defined a Silver \(d\)-cube \(Q_{m,n} = \prod\limits_{m\in d} [h_{m,n}]\)
  and numbers \( i_{m,0}, \ldots, i_{m,n} \not\in \dom(h_{m,n}) \).
  Fix a cover \( \mathcal{U}_n \) of \( X \) consisting of open sets of the diameter less than \( \frac{1}{n+1} \).
  As every metric space is paracompact, the cover \( \mathcal{U}_n \) has a \(\sigma\)-discrete refinement \( \tilde{\mathcal{U}}_n \).
  By the corollary above and \(\sigma\)-discreteness of \( \tilde{\mathcal{U}}_n \),
  for any Silver \(d\)-cube \( R \) there has to exist \( U \in \tilde{\mathcal{U}}_n \), such that \( U \cap R \in v^d \setminus v^d_0 \).
  Let \( (2^n)^m = \{ \alpha_k \colon k \in 2^{n \cdot m} \} \) and \( Q^{0}_{m,n} = [h^{0}_{m, n}] = Q_{m,n} \).
  For any \( k \in 2^n \) let
  \[
    h^{*k+1}_{m, n} = h^k_{m, n}|_{\omega \setminus I_{m,n})} \cup \{ (i_{m,0}, \alpha_{k}(m,0)), \ldots, (i_{m,n}, \alpha_{k}(m,n)) \},
  \]
  where \( I_{m,n} = \{ i_{m,0}, \ldots, i_{m,n} \} \).
  In a Silver \(d\)-cube \( \prod\limits_{m \in d} [h^{*k+1}_{m, n}] \)
  pick a Silver \(d\)-cube \( \prod\limits_{m \in d} [h^{k+1}_{m, n}] \)
  contained in some \(U \in \tilde{\mathcal{U_n}}\).
  We can take \( Q_{m, n+1} = \prod\limits_{m \in d} [h_{m, n+1}] = [h^{2^{n \cdot d}}_{m, n}|_{\omega \setminus I_{m,n}}] \)
  and \( i_{m,n+1} = \min (\omega \setminus D_{m,n+1}) \),
  where \( D_{m,n+1} = \dom (h_{m,n+1}) \).
  Clearly \( i_{m,n+1} > i_{m,n} \).
  As each set \( I_m = \bigcup\limits_{n \in \omega} I_{m,n} \) is infinite, we obtain that the set
  \[
    Q = \prod\limits_{m \in d} \bigcap\limits_{n \in \omega} Q_{m,n} = \prod\limits_{m \in d} [\bigcup\limits_{n \in \omega} h_{m,n}] =
    \prod\limits_{m \in d} [g_{m,n}]
  \]
  is a Silver \(d\)-cube.
  Furthermore, the sets of the form
  \begin{align*}
    Q \cap & \prod\limits_{m \in d} [h_{m, n} \cup \{ (i_{m,0}, \alpha_{k}(m,0)), \ldots, (i_{m,n}, \alpha_{k}(m,n)) \}] = \\
           & \prod\limits_{m \in d} [g_{m, n} \cup \{ (i_{m,0}, \alpha_{k}(m,0)), \ldots, (i_{m,n}, \alpha_{k}(m,n)) \}]
  \end{align*}
  form the base of topology on Q.
  It follows that the function \(f|_Q\) is continuous.

  For the case \(d = \omega\) let \(Q_{m,0}\) and \(i_{0,0} \min(\omega \setminus \dom(h_0))\).
  Suppose now that we have defined a Silver cube \(Q_{m,n} = \prod\limits_{m\in d} [h_{m,n}]\)
  and numbers \( i_{m,0}, \ldots, i_{m,n} \not\in \dom(h_{m,n}) \) for \(m \le n\).
  Once again, fix a cover \( \mathcal{U}_n \) of \( X \) consisting of open sets of the diameter less than \( \frac{1}{n+1} \),
  and take its \(\sigma\)-discrete refinement \( \tilde{\mathcal{U}}_n \).
  By the corollary above and \(\sigma\)-discreteness of \( \tilde{\mathcal{U}}_n \),
  for any Silver cube \( R \) there has to exist \( U \in \tilde{\mathcal{U}}_n \), such that \( U \cap R \in v^\omega \setminus v^\omega_0 \).
  Let \( (2^n)^n = \{ \alpha_k \colon k \in 2^{n^{2}} \} \) and \( Q^{0}_{m,n} = [h^{0}_{m, n}] = Q_{m,n} \).
  For any \( k \in 2^n \) let
  \[
    h^{*k+1}_{m, n} = h^k_{m, n}|_{(\omega \setminus I_{m,n})}
    \cup \{ (i_{m,0}, \alpha_{k}(m,0)), \ldots, (i_{m,n}, \alpha_{k}(m,n)) \}
  \]
  for \( m \le n\), and \(h^{*k+1}_{m, n} = h^{k}_{m, n}\) for \( m > n \),
  where \( I_{m,n} = \{ i_{m,0}, \ldots, i_{m,n} \} \).
  In a Silver cube \( \prod\limits_{m \in \omega} [h^{*k+1}_{m, n}] \)
  pick a Silver cube \( \prod\limits_{m \in \omega} [h^{k+1}_{m, n}] \)
  contained in some \(U \in \tilde{\mathcal{U_n}}\).
  We can take \( Q_{m, n+1} = \prod\limits_{m \in \omega} [h_{m, n+1}] = [h^{2^{n \cdot m}}_{m, n}|_{\omega \setminus I_{m,n}}] \)
  and \( i_{m,n+1} = \min (\omega \setminus D_{m,n+1}) \),
  where \( D_{m,n+1} = \dom (h_{m,n+1}) \).
  Clearly \( i_{m,n+1} > i_{m,n} \).
  Furthermore, take \( i_{0,n+1}, \ldots, i_{n+1, n+1} \) to the first \( n+1 \) natural numbers not in \( D_{n+1,n+1} \).
  As each set \( I_m = \bigcup\limits_{n \ge m} I_{m,n} \) is infinite and disjoint from any \( \dom(h_{m,n+1}) \), we obtain that the set
  \[
    Q = \prod\limits_{m \in \omega} \bigcap\limits_{n \in \omega} Q_{m,n} = \prod\limits_{m \in \omega} [\bigcup\limits_{n \in \omega} h_{m,n}]
  \]
  is a Silver cube.
  Moreover, the sets of the form
  \begin{align*}
    Q \cap & \prod\limits_{m \le n} [h_{m, n} \cup \{ (i_{m,0}, \alpha_{k}(m,0)), \ldots, (i_{m,n}, \alpha_{k}(m,n)) \}]
    \times \prod\limits_{m > n} [h_{m, n}]
    =                                                                                                                    \\
           & \prod\limits_{m \le n} [g_{m, n} \cup \{ (i_{m,0}, \alpha_{k}(m,0)), \ldots, (i_{m,n}, \alpha_{k}(m,n)) \}]
    \times \prod\limits_{m > n} [g_{m, n}]
  \end{align*}
  form the base of the product topology on Q.
  It follows that the function \(f|_Q\) is continuous.
\end{Prf}

A similar reasoning works for the ideal \( r^d_0 \) and the algebra \( r^d \).

\begin{Lem}
  Let \( \mathcal{F} = \{ F_i \colon i \in I \} \subseteq r^d_0 \) be a disjoint,  \( (r^d) \)-additive family. If there exist an Ellentuck \(d\)-cube \( V = \prod\limits_{m \in d} (s_m, A_m)^\omega \subseteq \bigcup \mathcal{F} \), then the set \(\{i \in I \colon F_i \cap V \not= \emptyset \}\) has cardinality \(\con \).
\end{Lem}
\begin{Prf}
  Without loss of generality, we can assume \( s_m = \emptyset \) for \( m \in \omega \).
  Similarly to perfect and Silver cubes, since \( r^d_0 \) is a \(\sigma\)-ideal in \( r^d \),
  the set \( \Delta = \{ J \subseteq I \colon \bigcup\limits_{i \in J} F_i \in r^d_0 \} \) is a free \(\sigma\)-ideal on \( I \).
  Let \(J = \{i \in I \colon F_i \cap V \not= \emptyset \}\).
  Clearly, \( J \not\in \Delta \) and there exist \( J_0, J_1 \subseteq J \) disjoint and both not in \( \Delta \).
  There exist Ellentuck \(d\)-cubes \( \prod\limits_{m \in d} (\emptyset, A^*_{m, 0})^\omega  \subseteq \bigcup\limits_{i \in J_0} F_i \)
  and \( \prod\limits_{m \in d} (\emptyset, A^*_{m, 1})^\omega \subseteq \bigcup\limits_{i \in J_1} F_i \).
  Let \( a_{m,0} = \min A^*_{m,0} \), \( A_{m,0} = A^*_{m,0} \setminus \{ a_{m,0} \} \), \( a_{m,1} = \min A^*_{m,1} \) and \( A_{m,1} = A^*_{m,1} \setminus \{ a_{m,1} \} \).
  Clearly, \( \prod\limits_{m \in d} (\{a_{m,0}\}, A_{m, 0})^\omega \subseteq \prod\limits_{m \in d} (\emptyset, A^*_{m, 0})^\omega \) and \( \prod\limits_{m \in d} (\{a_{m,1}\}, A_{m, 1})^\omega \subseteq \prod\limits_{m \in d} (\emptyset, A^*_{m, 1})^\omega \).

  Now assume that we have defined disjoint sets \( J_t \not\in \Delta \) for \( t \in 2^{n+1} \) as well as Ellentuck \(d\)-cubes
  \[
    \prod\limits_{m \in d} (\{ a_{m, t|1}, \ldots, a_{m, t|_{n+1}} \}, A_{m, t})^\omega  \subseteq \bigcup\limits_{i \in J_0} F_i.
  \]
  In each \( J_t \) we can find disjoint subsets \( J_{t^\smallfrown 0}, J_{t^\smallfrown 1} \),
  such that there exist Ellentuk \(d\)-cubes
  \( \prod\limits_{m \in d} (\{ a_{m, t|1}, \ldots, a_{m, t|_{n+1}} \}, A^*_{m, t^\smallfrown 0})^\omega \subseteq \bigcup\limits_{i \in J_{t^\smallfrown 0}} F_i \) and
  \( \prod\limits_{m \in d} (\{ a_{m, t|1}, \ldots, a_{m, t|_{n+1}} \}, A^*_{m, t^\smallfrown 1})^\omega \subseteq \bigcup\limits_{i \in J_{t^\smallfrown 1}} F_i \).
  We can take \( a_{m, t^\smallfrown 0} = \min A^*_{m, t^\smallfrown 0}\), \( A_{m, t^\smallfrown 0} = A^*_{m, t^\smallfrown 0} \setminus \{ a_{m, t^\smallfrown 0} \} \),
  \( a_{m, t^\smallfrown 1} = \min A^*_{m, t^\smallfrown 1}\) and \( A_{m, t^\smallfrown 1} = A^*_{m, t^\smallfrown 1} \setminus \{ a_{m, t^\smallfrown 1} \} \).

  For any \( x \in 2^\omega \) we have
  \[
    |\bigcap\limits_{n \in \omega} \prod\limits_{m \in d} (\{ a_{m, t|1}, \ldots, a_{m, t|_{n}} \}, A_{m, t|_{n}})^\omega| = 1.
  \]
  Therefore, \( \bigcap\limits_{n \in \omega} J_{t|_{n}} \not= \emptyset \), and thus \( |J| = \con \).
\end{Prf}

Just as in the case of \((s^d)\) and \((v^d)\) measurability applying a variant of the Berstein construction yields:
\begin{Cor}
  Let \( \mathcal{F} = \{ F_i \colon i \in I \} \subseteq r^d_0 \) be a disjoint,  \( (r^d) \)-additive family. Then \( \bigcup\mathcal{F} \in r^d_0 \).
\end{Cor}

\begin{Prf} of Theorem \ref{thm:r-property}.
  Let \( \prod_{m \in d} (s_m, A_m)^\omega \) be any Ellentuck \(d\)-cube.
  Without loss of generality, we can assume \( s_m = \emptyset \).
  Put \( A_{m,0} = A_m \).

  First consider the case \( d \in \omega \).
  Suppose we have defined the sets \( A_{m,n} \) as~well as~numbers~\( a_{m,k} \) for \( m \in d \) and \( k \in n \) such that
  \( \max(\{ a_{m,0}, \ldots, a_{m, n-1} \}) < \min(A_{m,n}) \).
  Like in the case of \((s^d)\) and \((v^d)\) measurable functions,
  fix a cover \( \mathcal{U}_n \) of \( X \) consisting of open sets of the diameter less than \( \frac{1}{n+1} \),
  and take its \(\sigma\)-discrete refinement \( \tilde{\mathcal{U}}_n \).
  By the corollary above and \(\sigma\)-discreteness of \( \tilde{\mathcal{U}}_n \),
  for any Ellentuck \(d\)-cube \( V \) there has to exist \( U \in \tilde{\mathcal{U}}_n \), such that \( U \cap V \in r^d \setminus r^d_0 \).
  Let \( \prod\limits_{m \in d}\Pow(\{ a_{m,0}, \ldots, a_{m, n-1} \}) = \{ \alpha^k \colon k \in 2^{d \cdot n} \} \)
  and \( A^0_{m,n} = A_{m,n} \).\\
  For~any~\(k \in d \cdot n\) in the Ellentuck \(d\)-cube \(\prod\limits_{m \in d} (\alpha^k_m, A^k_{m,n})^\omega\)
  find  the Ellentuck \(d\)-cube \(\prod\limits_{m \in d} (\alpha^k_m, A^{k+1}_{m,n})^\omega\)
  contained in some \(U \in \tilde{\mathcal{U_n}}\).
  Let \( a_{m,n} = \min A^{d \cdot n}_{m,n} \).
  We can take \( A_{m,n+1} = A^{2^{d \cdot n}}_{m,n} \setminus \{ a_{m,n} \} \).

  Let
  \[
    W = \prod\limits_{m \in d} (\emptyset, \{ a_{m, n} \colon n \in \omega \})^\omega = \prod\limits_{m \in d} (\emptyset, B_m)^\omega.
  \]
  Since the sets of the form
  \[
    W \cap \prod\limits_{m \in d} (\alpha^k_m, A{m,n+1})^\omega =
    \prod\limits_{m \in d} (\alpha^k_m, B_m)^\omega
  \]
  form the basis for the euclidean topology on \(B\),
  it follows that \( f|_B \) is continuous with respect to the euclidean topology.

  For the case \( d = \omega \) suppose we have defined the sets \( A_{m,n} \) for \( m \in \omega \) and numbers \( a_{m,k} \) for \( m, k \in n \), such that
  \( \max(\{ a_{m,0}, \ldots, a_{m, n-1} \}) < \min(A_{m,n}) \).
  Once again, fix a cover \( \mathcal{U}_n \) of \( X \) consisting of open sets having diameter less than \( \frac{1}{n+1} \),
  and take its \(\sigma\)-discrete refinement \( \tilde{\mathcal{U}}_n \).
  For any Ellentuck cube \( V \) there has to exist \( U \in \tilde{\mathcal{U}}_n \), such that \( U \cap V \in r^d \setminus r^d_0 \).
  Let \( \prod\limits_{m \in n}\Pow(\{ a_{m,0}, \ldots, a_{m, n-1} \}) = \{ \alpha^k \colon k \in 2^{n^2} \} \) and \( A^0_{m,n} = A_{m,n} \).
  For any \(k \in n^2\) in the Ellentuck cube \(\prod\limits_{m \in n} (\alpha^k_m, A^k_{m,n})^\omega \times \prod\limits_{m \ge n} (\emptyset, A^k_{m,n})^\omega \)
  find  the Ellentuck cube \(\prod\limits_{m \in n} (\alpha^k_m, A^{k+1}_{m,n})^\omega \times \prod\limits_{m \ge n} (\emptyset, A^{k+1}_{m,n})^\omega \)
  contained in some \(U \in \tilde{\mathcal{U_n}}\).
  Let \( a_{n, 0}, \ldots, a_{n,n} \) be the first \( n + 1 \) elements of \( A^{n^2}_{n,n} \),
  and for \( m \in n \) let \( a_{m,n} = \min A^{n^2}_{m,n} \).
  We can take \( A_{m,n+1} = A^{2^{n^2}}_{m,n} \setminus \{ a_{m,n} \} \) for \( m \in n \),
  \( A_{n, n+1} = A^{2^{n^2}}_{n,n} \setminus \{ a_{n, 0}, \ldots, a_{n,n} \} \),
  and \( A_{m, n+1} = A^{2^{n^2}}_{m,n} \) for \( m > n \).

  Let
  \[
    W = \prod\limits_{m \in \omega} (\emptyset, \{ a_{m, n} \colon n \in \omega \})^\omega
    = \prod\limits_{m \in \omega} (\emptyset, B_m)^\omega.
  \]
  Since the sets of the form
  \[
    W \cap \prod\limits_{m \in n} (\alpha^k_m, A{m,n+1})^\omega \times \prod\limits_{m \ge n} (\emptyset, A_{m,n+1})^\omega =
    \prod\limits_{m \in n} (\alpha^k_m, B_m)^\omega \times \prod\limits_{m \ge n} (\emptyset, B_{m})^\omega
  \]
  form the basis for the euclidean topology on \(B\),
  it follows that \( f|_B \) is continuous with respect to the euclidean topology.
\end{Prf}

%%%%%%%%%%%%%%%%%%%%
% 
% Application to Laver
% 
%%%%%%%%%%%%%%%%%%%%
\section{Application to generalization of Laver's theorem}

In~\cite{Lav}, Laver proved the following.
\begin{Thm}
  Let \( f_n \colon \prod\limits_{m \in \omega} Q_m \to [0;1] \) be all either continuous, Baire, or measurable function for \( n \in \omega \), where \( Q_m \) are perfect.
  Then, there exist a set \( N \in [\omega]^\omega \) as well as perfect sets \( P_m \subseteq Q_m \) for \( m \in \omega \),
  such that the sequence \( (f_n)_{n \in N} \) is monotonically (and thus uniformly) convergent on the product \( \prod\limits_{m \in d} P_m \).
\end{Thm}
It gave the positive answer to the question asked by Harrington in~\cite{AS}. Our results allow us to generalize this result to a wider class of functions.

\begin{Thm}
  Let \( f_n \colon \prod\limits_{m \in \omega} X_m \to X \) be \((s^\omega)\)-measurable functions.
  Then, there exists a perfect cube \( P = \prod\limits_{m \in \omega} P_m \),
  such that each \( f_n|_{P} \) is continuous.
\end{Thm}
\begin{Prf}
  Let the sets \( R_{m,0} \subseteq X_m \) be such
  that the function \( f_0 \) is continuous on the product \(\prod\limits_{m \in \omega} R_{m,0}\).
  Let \( P_{0, (0)} \) and \( P_{0, (1)} \) be two disjoint relative base sets in \( R_{0,0} \) of length at least \( 1 \).

  Assume inductively that for some \( n \in \omega \)
  we have defined the sets \( R_{m,n} \) for \( m > n \),
  and \( P_{m,t} \) for \( m \le n \) and \( t \in 2^{n+1} \),
  such that the functions \( f_k \) for \( k \le n \) are continuous on the set
  \[
    \prod\limits_{m \le n}\bigcup\limits_{t \in 2^{n+1}}P{m, t} \times \prod\limits_{m > n} R_{m,n}.
  \]
  Since there is finitely many sets \( P_{0, t_0} \times \ldots \times P_{n, t_n} \times \prod\limits_{m > n} R_{m,n} \),
  we can choose the sets \( R_{m,n+1} \subseteq R_{m,n} \) for \( m > n \),
  as well as \( Q_{m, t} \subseteq P_{m,t} \) for \( m \le n \) and \( t \in 2^{n+1} \),
  such that the function \( f_{n+1} \) is continuous on the set
  \[
    \prod\limits_{m \le n}\bigcup\limits_{t \in 2^{n+1}}Q_{m, t} \times \prod\limits_{m > n} R_{m,n+1}.
  \]
  In each set \( Q_{m, t} \) we can find two disjoint relative base sets \(P_{m, t^\smallfrown 0}, P_{m, t^\smallfrown 1}\) of length at least \( n + 1 \),
  and in the set \( R_{n+1,n+1} \) we can find \( 2^{n+2} \) disjoint relative base sets \(P_{n+1, t}\) of length at least \( n + 1 \).
  It follows that the functions \( f_n \) for \( k \le n+1 \) are continuous on the set
  \[
    \prod\limits_{m \le n + 1}\bigcup\limits_{t \in 2^{n+2}}P_{m, t} \times \prod\limits_{m > n + 1} R_{m,n+1}.
  \]

  Now let
  \[
    P_m = \bigcap\limits_{n \ge m} \bigcup\limits_{t \in 2^{n + 1}} P_{n, t}.
  \]
  Clearly each set \( P_m \) is perfect, and all the functions \( f_n \) are continuous on the cube \(\prod\limits_{m \in \omega}P_m\).
\end{Prf}

\begin{Thm}
  Let \( f_n \colon (2^\omega)^\omega \to X \) be \((v^\omega)\)-measurable functions.
  Then, there exists a Silver cube \( Q = \prod\limits_{m \in \omega} Q_m \),
  such that each \( f_n|_{P} \) is continuous.
\end{Thm}
\begin{Prf}
  Let \(Q_{m,0} = 2^\omega\) and \(i_{0,0} = 0\).
  Suppose now that we have already defined a Silver cube \(Q_{n} = \prod\limits_{m\in d} [h_{m,n}]\)
  and numbers \( i_{m,0}, \ldots, i_{m,n} \not\in \dom(h_{m,n}) \) for \(m \le n\),
  such that \(f_k|_{Q_n}\) is continuous for \( k \in n \).
  Let \( (2^n)^n = \{ \alpha_k \colon k \in 2^{n^{2}} \} \) and \( Q^{0}_{m,n} = [h^{0}_{m, n}] = Q_{m,n} \).
  For any \( k \in 2^n \) let
  \[
    h^{*k+1}_{m, n} = h^k_{m, n}|_{(\omega \setminus I_{m,n})}
    \cup \{ (i_{m,0}, \alpha_{k}(m,0)), \ldots, (i_{m,n}, \alpha_{k}(m,n)) \}
  \]
  for \( m \le n\), and \(h^{*k+1}_{m, n} = h^{k}_{m, n}\) for \( m > n \),
  where \( I_{m,n} = \{ i_{m,0}, \ldots, i_{m,n} \} \).
  In a Silver cube \( \prod\limits_{m \in \omega} [h^{*k+1}_{m, n}] \)
  pick a Silver cube \( \prod\limits_{m \in \omega} [h^{k+1}_{m, n}] \),
  such that \(f_n\) is continuous on it.
  We can take \( Q_{m, n+1} = \prod\limits_{m \in \omega} [h_{m, n+1}] = [h^{2^{n \cdot m}}_{m, n}|_{\omega \setminus I_{m,n}}] \)
  and \( i_{m,n+1} = \min (\omega \setminus D_{m,n+1}) \),
  where \( D_{m,n+1} = \dom (h_{m,n+1}) \).
  Clearly \( i_{m,n+1} > i_{m,n} \).
  Furthermore, take \( i_{0,n+1}, \ldots, i_{n+1, n+1} \) to the first \( n+1 \) natural numbers not in \( D_{n+1,n+1} \).
  As each set \( I_m = \bigcup\limits_{n \ge m} I_{m,n} \) is infinite and disjoint from any \( \dom(h_{m,n+1}) \) we obtain that the set
  \[
    Q = \prod\limits_{m \in \omega} \bigcap\limits_{n \in \omega} Q_{m,n} = \prod\limits_{m \in \omega} [\bigcup\limits_{n \in \omega} h_{m,n}]
  \]
  is a Silver cube.
  Moreover, each of the functions \(f_n|_Q\) is continuous.
\end{Prf}

\begin{Thm}
  Let \( f_n \colon ([\omega]^{\omega})^\omega \to [0;1] \) be \((r)^\omega\)-measurable functions.
  Then, there exists an Ellentuck cube \( \prod_{m \in \omega} (\emptyset, A_m)^\omega \),
  such that each \( f_n|_{\prod_{m \in \omega} (\emptyset, A_m)^\omega} \) is continuous with respect to the euclidean topology.
\end{Thm}
\begin{Prf}
  Let \(A_{m,0} = \omega\).
  Suppose we have defined the sets \( A_{m,n} \) for \( m \in \omega \), and numbers \( a_{m,k} \) for \( m, k \in n \) such that,
  \( \max(\{ a_{m,0}, \ldots, a_{m, n-1} \}) < \min(A_{m,n}) \).

  Let~\( \prod\limits_{m \in n}\Pow(\{ a_{m,0}, \ldots, a_{m, n-1} \}) = \{ \alpha^k \colon k \in 2^{n^2} \} \) and \( A^0_{m,n} = A_{m,n} \).
  For any \(k \in n^2\) in the Ellentuck cube
  \[\prod\limits_{m \in n} (\alpha^k_m, A^k_{m,n})^\omega \times \prod\limits_{m \ge n} (\emptyset, A^k_{m,n})^\omega \]
  find  the Ellentuck cube
  \[\prod\limits_{m \in n} (\alpha^k_m, A^{k+1}_{m,n})^\omega \times \prod\limits_{m \ge n} (\emptyset, A^{k+1}_{m,n})^\omega \]
  on which the function \(f_n\) is continuous.
  Let \( a_{n, 0}, \ldots, a_{n,n} \) be the first \( n + 1 \) elements of \( A^{n^2}_{n,n} \),
  and \( m \in n \) let \( a_{m,n} = \min A^{n^2}_{m,n} \).
  We can take \( A_{m,n+1} = A^{2^{n^2}}_{m,n} \setminus \{ a_{m,n} \} \) for \( m \in n \),
  \( A_{n, n+1} = A^{2^{n^2}}_{n,n} \setminus \{ a_{n, 0}, \ldots, a_{n,n} \} \),
  and \( A_{m, n+1} = A^{2^{n^2}}_{m,n} \) for \( m > n \).

  Let
  \[
    W = \prod\limits_{m \in \omega} (\emptyset, \{ a_{m, n} \colon n \in \omega \})^\omega
    = \prod\limits_{m \in \omega} (\emptyset, A_m)^\omega.
  \]
  We obtain that all the functions \( f_n \) are continuous with respect to the euclidean topology on \( W \).
\end{Prf}

\begin{Cor}
  Let \( f_n \colon [\omega]^{\omega} \to [0;1] \) be \(CR\)-measurable functions. Then there exists a set \( P \subseteq [\omega]^{\omega} \) homeomorphic to \( 2^\omega \), such that each \( f_n|_{P} \) is continuous with respect to the euclidean topology.
\end{Cor}
\begin{Prf}
  Each set \( (s, A)^\omega \) is homeomorphic in euclidean topology to the space of irrational numbers \( \omega^\omega \). The result follows in a straightforward way.
\end{Prf}

\end{document}